
%
\magnification=\magstep1
\newif\ifmssymb  %
\mssymbtrue      
\input mssymb    
\catcode`\@=11   %
\def\newfam{\alloc@8\fam\chardef\sixt@@n}
\catcode`\@=12
\def\today{November 15, 1993}
\font\petcap=cmcsc10
\font\sevenit=cmti7 scaled \magstep0
\font\sevensl=cmti7 scaled \magstep0
\font\bigit=cmti10 scaled \magstep1
\font\bigbf=cmbx10 scaled \magstep1
\font\bigrm=cmr10 scaled \magstep1
\font\bigsevenbf=cmbx7 scaled \magstep1
\font\bigsevenrm=cmr7 scaled \magstep1
\font\bigteni=cmmi10 scaled \magstep 1
\font\bigex=cmex10 scaled \magstep1
\font\bigsy=cmsy10 scaled \magstep1
\font\bigseveni=cmmi7 scaled \magstep1
\font\bigsevensy=cmsy7 scaled \magstep1
\ifmssymb
\font\bigmsy=msbm10 scaled \magstep1
\font\bigsevenmsy=msbm7 scaled \magstep1
\font\bigmsx=msam10 scaled  \magstep1
\fi
\skewchar\bigteni='177
\skewchar\bigsy='60
\catcode`\@=11
\newskip\standardbaselineskip
\standardbaselineskip=14pt
\def\tenpoint{%
 \textfont0=\tenrm
 \textfont1=\teni
 \textfont2=\tensy
 \textfont\bffam=\tenbf
 \textfont\ttfam=\tentt
 \textfont\itfam=\tenit
 \textfont\slfam=\tensl
\ifmssymb
\textfont\msxfam=\tenmsx
\textfont\msyfam=\tenmsy
\fi
\scriptfont0=\sevenrm \scriptscriptfont0=\fiverm
 \scriptfont1=\seveni \scriptscriptfont1=\fivei
 \scriptfont2=\sevensy \scriptscriptfont2=\fivesy
 \scriptfont3=\tenex \scriptscriptfont3=\tenex
 \scriptfont\bffam=\sevenbf \scriptscriptfont\bffam=\fivebf
\scriptfont4=\sevenit
\scriptfont5=\sevensl
 \ifmssymb\scriptfont\msyfam=\sevenmsy \scriptscriptfont\msyfam=\fivemsy
 \scriptfont\msxfam=\sevenmsx \scriptscriptfont\msxfam=\fivemsx
 \fi
  \def\it{\fam\itfam\tenit}%
  \def\rm{\fam\z@\tenrm}%
  \def\bf{\fam\bffam\tenbf}%
  \def\sl{\fam\slfam\tensl}%
  \def\tt{\fam\ttfam\tentt}%
  \def\oldstyle{\fam\@ne\teni}%
  \def\big##1{{\hbox{$\left##1\vbox to 8.5pt{}\right.\n@space$}}}
  \setbox\strutbox=\hbox{\vrule height8.5pt depth3.5pt width0pt}%
  \abovedisplayskip=13pt plus 3pt minus 10pt
  \belowdisplayskip=13pt plus 3pt minus 10pt
  \normalbaselineskip=\standardbaselineskip
  \normalbaselines\rm}
\def\ninepoint{\tenpoint}
\def\bigstyle{
\textfont0=\bigrm
\scriptfont0=\bigsevenrm
\textfont1=\bigteni
\scriptfont1=\bigseveni
\textfont2=\bigsy
\scriptfont2=\bigsevensy
\textfont3=\bigex
\textfont4=\bigit
\scriptfont4=\tenit
\textfont5=\bigit %
\textfont6=\bigbf
\scriptfont6=\bigsevenbf
\textfont\ttfam=\nullfont
\ifmssymb
\textfont\msyfam=\bigmsy
\scriptfont\msyfam=\bigsevenmsy
\textfont\msxfam=\bigmsx
\fi
  \def\it{\fam\itfam\bigit}%
  \def\rm{\fam\z@\bigrm}%
  \def\bf{\fam\bffam\bigbf}%
  \def\sl{\fam\slfam\bigit}%
  \def\tt{\errmessage{\string\tt\ in \string\bigstyle\ not yet implemented}}%
  \def\oldstyle{\fam\@ne\bigteni}%
  \def\big##1{{\hbox{$\left##1\vbox to 10.2pt{}\right.\n@space$}}}
  \setbox\strutbox=\hbox{\vrule height10.2pt depth4.2pt width0pt}%
  \normalbaselineskip=\standardbaselineskip
  \advance\normalbaselineskip by 3pt
  \normalbaselines\rm}
\catcode`\@=12
\tenpoint
\ifmssymb
\def\Bbb#1{{\tenmsy\fam\msyfam\relax #1}}\else
\def\Bbb#1{{\bf #1}}
\fi
\def\goth#1{{\bf #1}}
\def\eg{e.g.\spacefactor=1000\ }%
\def\ie{i.e.\spacefactor=1000\ }%
\let\oldcal\cal
\def\cal#1{{\oldcal #1}}%
\let\Cedille\c
\def\a#1{{\bf #1}}%
\def\abs#1{{\vert #1\vert}}

\def\dimc{\mathop{{\fam=-1 dim}_{\s C}}}

\def\sdr{\mathbin{\vrule height 4.5pt depth 0.1pt\kern-1.8pt\times}}
\def\={\hbox{--}\nobreak\hskip0pt\relax}
\def\({\left(}
\def\){\right)}%
\pageno=1
\relpenalty=5000
\binoppenalty=5000
\hfuzz=2truept
\mathsurround=1pt
\nulldelimiterspace=0pt
\overfullrule=0pt
\normallineskip=1pt
\normallineskiplimit= 1pt
\normalbaselines
\newdimen\oldparindent

\def\restoreparindent{\let\parindent\oldparindent}%
\def\.{\null.}%
\newtoks\autor
\autor={J\"org Winkelmann}%
\newtoks\ort
\ort={Bochum}%
\newcount\arbvar
\newtoks\normalheadline
\normalheadline={\hfil}
\headline={\ifnum\pageno=1\hfill\else\the\normalheadline\fi}%
\catcode`\@=11
\def\makeheadline{\vbox to\z@{\vskip-32.5\p@
  \line{\vbox to8.5\p@{}\the\headline}\vss}\nointerlineskip}%
\catcode`\@=12
\newtoks\normalfootline
\footline={\hfill\tenrm\folio\hfill}
\def\finalversion{\relax}
\newskip\vardimen
\newskip\varlastdimen
\def\varskip{\varvarskip{0pt}}%
\def\varvarskip#1#2\penalty#3 {\varlastdimen=\lastskip
\removelastskip\penalty#3\relax
\vardimen=#2\relax
\advance\varlastdimen by#1\relax
\ifdim\vardimen<\varlastdimen\relax\vardimen=\varlastdimen\fi
\vskip\vardimen\relax}%
\catcode`\@=11
\def\footnote#1{\let\@sf\empty %
  \ifhmode\edef\@sf{\spacefactor\the\spacefactor}\/\fi
 $^{\rm #1}$\@sf\vfootnote{#1}}
\def\vfootnote#1{\insert\footins\bgroup
  \interlinepenalty\interfootnotelinepenalty
  \ninepoint
  \splittopskip\ht\strutbox %
  \splitmaxdepth\dp\strutbox \floatingpenalty\@MM %
  \leftskip\z@skip \rightskip\z@skip \spaceskip\z@skip \xspaceskip\z@skip
  \textindent{$^{\rm #1}$}\footstrut\futurelet\next\fo@t}%
\catcode`\@=12
\newtoks\postdisplaytoks
\everydisplay=\expandafter{\the\everydisplay\global\postdisplaytoks={}%
\aftergroup\thepostdisplaytoks}%
\def\thepostdisplaytoks{\the\postdisplaytoks}%
\def\afterdisplay#1{\global\postdisplaytoks=\expandafter{\the\postdisplaytoks
#1}}%
\newif\ifdisplay
\global\displayfalse
\everydisplay=\expandafter{\the\everydisplay\displaytrue}%
\newcount\qedpenalty
\def\qed{\qedsign\ifdisplay\postdisplaypenalty=\qedpenalty
\afterdisplay\endproofstyle
\let\next\relax\else\par\penalty\qedpenalty\let\next\endproofstyle\fi\next}%
\def\qedsign{\ifmmode\ \qedkasten{4pt}\else\nobreak\ $\qedkasten{4pt}$\fi}%
\def\qedkasten#1{\mathop{\kern0.5pt\vbox{\hrule\hbox{\vrule%
   \hskip#1\vrule height #1 width0pt%
   \vrule}\hrule}\kern-0.5pt}}%
\qedpenalty=-100
\def\endproofstyle{\tenpoint\medskip}%
\ifmssymb

\fi
\def\a#1{{\goth #1}}%
\def\s#1{\Bbb #1}%
\normalheadline=%
{\ninepoint\sl\hfil\ifodd\pageno\the\shorttitle
\else\the\autor\fi\hfil}%
\newbox\titbox
\def\fintitre{\global\let\titre\relax}%
\def\titel#1{\bigskip\bigskip
\def\titre##1\\{
\setbox\titbox=%
\centerline{\bigstyle\bf ##1}%
\ifdim\ht\titbox>0pt\relax \box\titbox\medskip\fi
\titre}%
\titre#1\fintitre\\
\medskip
\centerline{\sl by}
\medskip
\centerline{\rm\the\autor}
\bigskip\bigskip
}%
\newtoks\shorttitle
\shorttitle={}%
\newtoks\autoraddress
\autoraddress={%
J\"org Winkelmann\par
Mathematisches Institut\par
NA 4/69 \par
Ruhr-Universit\"at Bochum \par
D-44780 Bochum \par
Germany \par
\medskip
e-mail (BITNET):\par
WINKELMANN\par
@RUBA.RZ.RUHR-UNI-BOCHUM.DE\par
}%
\def\datum{\par\medskip\bgroup\leftskip=200pt \parindent=0pt \the\ort,
\today\par\bigskip\egroup}
\def\endarticle{\penalty100\relax\bigskip\datum\penalty-300
\smallskip %
\theaddress\end}%
\def\theaddress{%
\vbox{\leftskip=200pt
\parindent=0pt %
\baselineskip=8.5pt
\interlinepenalty=4000
\rm\the\autoraddress}}
\newtoks\thankstoks
\def\thanks#1{\thankstoks={#1}}
\def\thethanks{\par\the\thankstoks\par}
\newcount\sectno
\sectno=0
\def\section#1{\varskip\bigskipamount\penalty-250
\advance\sectno by 1
\leftline{\ifnum\sectno>0\relax\the\sectno. \else\fi
\petcap\ignorespaces #1}%
\nobreak
\medskip}%
\def\subsection#1{\varskip\medskipamount\penalty-120
\centerline{\petcap #1}\nobreak\medskip}%
\def\subsubsection#1{\varskip\medskipamount\penalty-120
\bgroup\it\ignorespaces#1. \egroup\ignorespaces}
\catcode`\@=11
\bgroup\obeylines\gdef\proclaim@@#1#2#3#4#5
{\egroup\varskip\bigskipamount\penalty30 %
\noindent{#1#2%
\def\parameter{#5}%
\ \ignorespaces\ifx\parameter\empty #3\else #5\fi\unskip. }%
\bgroup#4\ignorespaces}\egroup
\def\endproclaim{\par\egroup\penalty-50\medskip}%
\def\proclaim@{\bgroup\obeylines\proclaim@@}%
\def\proclaim#1{\proclaim@{\bf}{#1}{}{\it}}
\def\Corollary{\proclaim@{\bf}{Corollary}{}{\it}}
\def\Example{\proclaim@{\bf}{Example}{}{\rm}}
\def\Definition{\proclaim@{\bf}{Definition}{}{\rm}}
\def\Remark{\proclaim@{\bf}{Remark}{}{\rm}}
\def\Proof{\proclaim@{\it}{Proof}{}{\ninepoint\rm}}
\newcount\propno
\propno=1
\def\Proposition{\proclaim@{\bf}{Proposition}
{\the\propno\global\advance\propno by 1 }{\it}}
\newcount\theono
\theono=1
\def\Theorem{\proclaim@{\bf}{Theorem}
{\the\theono\global\advance\theono by 1 }{\it}}
\newcount\lemno
\lemno=1
\def\Lemma{\proclaim@{\bf}{Lemma}{\the\lemno\global\advance\lemno by 1 }{\it}}
\catcode`\@=12

\def\endremark{\endproclaim}%
\def\endproofstyle{\egroup}%
\newif\ifciteerror
\newcount\citeno
\def\checkciteerror{\ifciteerror\errmessage{Citation error}\fi}
\def\defcitations#1{\dodefcites#1,,\enddefcite}
\def\dodefcites#1,{\def\next{#1}\ifx\next\empty
\let\next\finishdefcite
\else\let\next\doodefcite\fi
\next{#1}}
\def\doodefcite#1{\advance\citeno by 1
\expandafter\xdef\csname cite:#1\endcsname
{\noexpand\global\noexpand\citeerrorfalse
\the\citeno}\dodefcites} \def\finishdefcite#1\enddefcite{}
\def\makecite#1#2{\csname cite:#1\endcsname#2}
\def\docite#1,#2,#3!{\def\next{#2}\ifx\next\empty[\makecite{#1}{}]\else
[\makecite{#1}{,#2}]\fi}
\def\cite#1{{\rm\citeerrortrue\docite#1,,!\checkciteerror}}
\newdimen\refindent
\def\References{\medskip\sectno=-100 \section{References}
\let\c\Cedille
\def\(##1\){}
\frenchspacing
\setbox0=\hbox{99.\enspace}\refindent=\wd0\relax
\def\item##1##2:{\goodbreak
\hangindent\refindent\hangafter=1\noindent\hbox
to\refindent{\hss\refcite{##1}\enspace}\petcap
\ignorespaces##2\rm: }}
\def\refcite#1{\citeerrortrue\dorefcite#1[]\endcite\checkciteerror.}
\def\dorefcite#1[#2]#3\endcite{\def\next{#2}%
\ifx\next\empty\makecite{#1}{}\else\makecite{#2}{}\fi}%
\def\AMS#1{\relax}
\AMS{14G20,14L10,20G25,22E40}
\def\tensor{\otimes}
\def\dosmallmatrix#1&#2\cr#3&#4\cr{\bigl({#1\atop#3 }{#2\atop#4}\bigr)}%
\def\implies{\Rightarrow}
\def\Lie{\mathop{{\cal L}{\sl ie}}}
\finalversion
\defcitations{B,BS,BT,G,Mc,M,OP,R,Ro,S,T,Z}
\titel{On discrete Zariski-dense subgroups of algebraic groups}
\section{Introduction}
Let $k$ be a local field, \ie a field equipped with an absolute value
inducing a non-discrete locally compact topology.
(It is well-known that any local field $k$ is isomorphic to one of
the following: $\s C$, $\s R$, a finite extension of a field $\s Q_p$
of p-adic numbers or $\s F_q((t))$, where $\s F_q$ is a finite field.)
Then every $k$-variety admits two topologies: the Zariski-topology
and a Hausdorff topology induced by the absolute value on $k$.
We are now interested in subgroups of $k$-groups for which these two
topologies are as different as possible, \ie groups which are discrete
in the Hausdorff topology, but dense in the Zariski-topology.

A main motivation for our investigation was the following.
For complex linear-algebraic groups a discrete cocompact subgroup
is necessarily Zariski-dense. There are known obstructions to the
existence of discrete cocompact subgroups, namely only unimodular
groups may admit discrete cocompact subgroups.
Thus one may ask whether these obstructions actually prohibit only
discrete cocompact subgroups or prohibit {\sl all} Zariski-dense
discrete subgroups. It turns out that there are many groups
(in particular all non-solvable complex parabolic groups) which
do admit discrete Zariski-dense subgroups although they do not
admit any discrete cocompact subgroup.
On the other hand there do exist linear-algebraic groups
without any discrete Zariski-dense subgroup at all.
The goal of this article is to give criteria as precise as possible
for the existence and non-existence of discrete Zariski-dense subgroups.

For non-solvable groups the main result is:
\Theorem
Let $G$ be a Zariski-connected non-solvable $k$-group,
defined over a local field $k$.
Let $R$ denote the radical.
For $char(k)=0$ the following statements are equivalent.
\item{(i)} $G/R(k)$ is non-compact in the Hausdorff topology.
\item{(ii)} $G(k)/R(k)$ is non-compact in the Hausdorff topology.
\item{(iii)} $G/R$ is $k$-isotropic, \ie contains a
positive-dimensional $k$-split torus.
\item{(iv)} $G(k)$ is not amenable.
\item{(v)} There exists a discrete Zariski-dense subgroup in $G(k)$.
\endproclaim

For $k=\s C$ this boils down to the following
\Theorem 1'
Let $G$ be a non-solvable connected complex linear-algebraic group.
Then $G$ admits a subgroup $\Gamma$ which is discrete in the Hausdorff
topology and dense in the Zariski-topology.
\endproclaim
Since parabolic complex groups never contain lattices
(they are never unimodular), this implies that
there exist complex
groups which admit discrete Zariski-dense subgroups although they
do not admit lattices.

For $char(k)>0$ we have partial results, in particular the implications
$ (v)\implies (iv)\implies (ii)$ of Theorem 1
still hold.

Concerning solvable groups we obtain quite different results for the
various local fields.
In the complex case, \ie $k=\s C$
we have a number of partial results.
\Theorem
\def\f{\par\noindent\sl}
{\f Reductive groups.}
Let $G$ be a connected reductive solvable group, \ie $(\s C^*)^n$.
Then for each $1\le k\le n $ there exists a Zariski-dense discrete
subgroup $\Gamma$ in $G$ with $\Gamma\simeq(\s Z^k,+)$.
(This is well-known). 
{\f Unipotent groups.}
A unipotent complex group $G$ admits a Zariski-dense discrete subgroup
if and only if $G$ admits a real-algebraic subgroup $G_0$ which may be
defined over $\s Q$ and is not contained in any proper complex subgroup
of $G$.
{\f Borel groups.}
A Borel group in a simple complex group never admits a Zariski-dense
discrete subgroup. In contrast, Borel groups in $SL_2(\s C)\times SL_2(\s C)$
contain Zariski-dense discrete subgroups.
{\f Metabelian groups.}
\item{(i)} Assume that $G=(\s C^*)\sdr_\rho(\s C^k,+)$.
Then unimodularity (\ie $\rho(\s C^*)\subset SL_k(\s C)$) is a necessary
but in general not sufficient condition for the existence of a Zariski-dense
discrete subgroup.
\item{(ii)}
A Borel group $B$ in $SL_2(\s C)\times SL_2(\s C)$ is an example for a
semidirect product $B\simeq(\s C^*)^2\sdr(\s C^2,+)$ which is not
unimodular but nevertheless admits a discrete Zariski-dense subgroup.
\item{(iii)} (Otte-Potters):
Let $T$ be a maximal torus in $S=SL_k(\s C)$ and $G=T\sdr_\rho(\s C^k,+)$
with $\rho$ given by the usual $S$-action on $\s C^k$.
Then $G$ admits a Zariski-dense discrete cocompact subgroup.
\endproclaim

For $k=\s R$ we obtain the following
\Theorem
Let $G$ be a connected solvable $\s R$-group.
Then the following conditions are necessary,
but in general not sufficient for the existence of
discrete Zariski-dense subgroups in $G(\s R)$.
\item{(i)} The commutator group $\cal D(G)$ may be defined over $\s Q$.
\item{(ii)} The group is unimodular, \ie $Ad(G)\subset SL(\a g)$.

If $G$ is unipotent, then condition $(i)$ is sufficient for the
existence of discrete Zariski-dense subgroups
(Condition $(ii)$ is automatically fulfilled, if $G$ is nilpotent).

\endproclaim
Finally, for non-archimedean fields
in characteristic zero we are able to provide a complete description.
\Theorem
Let $G$ be a Zariski-connected solvable $k$-group,
$k$ a non-archi\-me\-dean local field with
$char(k)=0$. Let $U$ denote the unipotent radical and $T$ a maximal
$k$-split torus.
Then there exists a discrete Zariski-dense subgroup of $G(k)$ if and only
if the following two conditions are fulfilled:
\item{(i)}
$G$ is commutative.
\item{(ii)}
$\dim T\ge\max\{1,\dim U\}$.
\endproclaim

We basically use the notation of \cite B and \cite M.
In particular a $k$-group $G$ is a linear-algebraic group defined over
a field $k$ and $G(k)$ denotes the group of $k$-rational points.
For a local field $k$, $G(k)$ carries an induced Hausdorff topology.
Topological terms refer to this Hausdorff topology.
Topological notions concerning the Zariski topology are preceded
by the prefix {\sl "Zariski-"}.
For $k=\s C$ the notions {\sl connected} and {\sl Zariski-connected}
coincide.
\section{Arithmetic Groups}
Let $k$ be a local field and $G$ a simple group defined over $k$.
Then there exists an {\sl arithmetic} subgroup $\Gamma\subset G(k)$,
which is a lattice (\ie discrete with finite covolume)
by the Borel-Harish-Chandra-Behr-Harder
reduction theorem. 
Due to the Borel-Wang Density theorem
$\Gamma$ is Zariski-dense unless $G$ is $k$-anisotropic.
(For a local field $k$ a reductive group $G$ is
$k$-anisotropic if and only if $G(k)$ is compact in the
Hausdorff topology \cite{M,2.3.6. p.54}).
Together this yields the following
\Theorem A
Let $G$ be a Zariski-connected
semisimple $k$-group without $k$-anisotropic factors,
$k$ a local field.
Then $G(k)$ contains a discrete Zariski-dense subgroup which furthermore
is a lattice.
\endproclaim
In particular every connected complex semisimple group and every
connected real semisimple
group without compact factor admit a discrete Zariski-dense subgroup.

\section{Preparations}

We need the following theorem of Tits on the existence of free subgroups
in linear groups.
\Theorem B (Tits)
Let $G$ be a semisimple
linear-algebraic group defined over a local field $k$,
$\Lambda$ a Zariski-dense subgroup of $G(k)$.
Then $\Lambda$ contains a free subgroup $H$ with
infinitely many generators $f_i$, such that any two of these $f_i$
generate a Zariski-dense subgroup of $G$.
\endproclaim
\Corollary
Let $G$ be a linear-algebraic group defined over a local field $k$,
$\Lambda$ a subgroup.
Then either $\Lambda$ contains a non-commutative free subgroup or it contains a
solvable subgroup of finite index.
\endproclaim
For $char(k)=0$ the theorem is Th.3 in \cite T.
For $char(k)>0$ it follows from Th. 4 in \cite T,
because in a local field of positive characteristic only finitely many
elements are algebraic over the prime field. (This follows from the result
that such a field is isomorphic to $\s F_q((t))$ for some finite field
$\s F_q$.)

Note that such a statement does not hold for arbitrary (\ie not local)
fields.
For instance, let $k$ be an algebraic closure of a finite field $\s F_p$.
Then $G=SL_2(k)$ is a group which contains neither a solvable subgroup
of finite index nor a free subgroup. In fact any finitely generated
subgroup of $G$ is finite. See \cite T for details.
We will need a result in order to control the image of
$k$-rational points under a $k$-morphism.
\Theorem C (Borel-Serre \cite{BS,p.153})
Let $k$ be a local field of characteristic zero, $G$, $H$ $k$-groups,
$\rho:G\to H$ a surjective $k$-group homomorphism.
Then $\rho(G(k))$ is a subgroup of finite index in $H(k)$.
\endproclaim
Unfortunately such a statement does not hold in positive
characteristic.
\Example 1
Let $k=\s F_p((t))$, $G_m$ the multiplicative group of the field
(\ie $G_m(k)=(k^*,\cdot)$) and $\rho:G_m\to G_m$ the group morphism
given by $\rho(x)=x^p$. Then $\rho(G_m(k))=\{\sum_k a_kt^{kp}\}$ and
the natural group homomorphism $\pi:G_m(k)\to G_m(k)/\rho(G_m(k))=Q$
maps $S=\{1+t+\sum_{k>0} a_kt^{kp}\}$ injectively into $Q$,
hence $\rho(G_m(k))$ has infinite index in $\rho(G_m)(k)$.
\endproclaim

\section{Non-amenable Groups}

\Theorem
Let $k$ be a local field of characteristic zero, $G$ be a
Zariski-connected $k$-group.
Assume that $G/R(k)$ is not compact, \ie not $k$-anisotropic.

Then $G(k)$ admits a Zariski-dense discrete subgroup.
\endproclaim
\Proof
By assumption there exists a surjective $k$-group morphism $\rho$
from $G$ to a simple $k$-isotropic group $S$.
$S(k)$ admits a Zariski-dense discrete subgroup $\Lambda$.
$I=\rho(G(k))$ is a subgroup of finite index in $S(k)$.
Thus $\Lambda_1=\Lambda\cap I$ is of finite index in $\Lambda$,
hence Zariski-dense and discrete in $S(k)$.
Now $\Lambda_1$ contains
a countable infinite subset $F=\{a_i:i\in\s N_0\}$
such that the elements of $F$ are
free generators of a discrete subgroup $\Lambda_0\subset\Lambda_1$ and
any two elements of $F$ generate a Zariski-dense subgroup of $S$.
Choose $b_0,b_1\in G(k)$ such that $\rho(b_i)=a_i$ ($i=1,2$) and let
$H$ denote the Zariski-closure of the subgroup generated by $b_0$ and
$b_1$. Now $\rho$ maps $H$ surjectively on $S$,
since $a_0$ and $a_1$
generate a Zariski-dense subgroup of $S$.
Hence $\rho(H(k))$ is a subgroup of finite index in $S(k)$.
Choose $c_i\in H(k)$, $n_i\in\s N$ ($i\ge 2$),  such that
$\rho(c_i)=a_i^{n_i}$.
Let $\Sigma=\{s_i:i\in\s N, i\ge 2\}$
be a countable Zariski-dense subset of $A(k)$ where
$A$ denotes the kernel of $\rho:G\to S$.
Finally choose $b_i\in G(k)$ for $i\ge 2$ such that $b_i=c_i\cdot s_i$
and let $\Gamma$ denote the subgroup of $G(k)$ generated by
the elements $b_i$ ($i\ge 0$). The map $\rho$ maps $\Gamma$ injectively
into $\Lambda$, hence $\Gamma$ is discrete.
Futhermore the construction implies that $\Gamma$ is Zariski-dense
in $G$.
\qed
In the above proof we used the assumption $char(k)=0$ at two points:
First we used that for a perfect field $k$ the radical of a $k$-group
is defined over $k$.
(A local field $k$ is perfect if and only if $char(k)=0$.)
Second we used that for a surjective $k$-group homomorphism
$\rho:G\to H$ the image of the $k$-rational points $\rho(G(k))$ has
finite index in $H(k)$.
By strengthening the assumptions one may circumvent these problems
and obtain the following result.
\Proposition
Let $k$ be a local field of arbitrary characteristic,
$S$ a $k$-isotropic simple $k$-group, $H$ a $k$-group and $G$ a
$k$-group which is $k$-isomorphic to a semi-direct product $S\sdr H$.

Then $G(k)$ admits a Zariski-dense discrete subgroup.
\endproclaim
\Proof
The assumptions imply the existence of a $k$-group homomorphism $\rho$ from
$G$ to $S$ which maps $G(k)$ surjectively on $S(k)$. Thus the proof
for $char(k)=0$ can be carried over with the following modification.
Using the semi-direct product
structure we may choose $b_0$, $b_1$ inside the factor isomorphic to $S$.
Then $H$ contains the factor isomorphic to $S$ and consequently $H(k)$ maps
surjectively on $S(k)$ by $\rho$.
\qed
\medskip 
\section{Amenable Groups}
\Definition
A locally compact topological group is called {\sl amenable} if for
every continuous action on a compact metrizable topological space there
exists an invariant probability measure.
\endproclaim
For any $k$-group defined over a local field $k$ the group $G(k)$ of
$k$-rational points is a locally compact topological group in the
Hausdorff topology. Thus we may apply the theory of amenable groups.
We summarize some important properties of amenable groups.
(see \cite G \cite Z for proofs and details).
\item{i)} Compact groups are amenable,
\item{ii)} Solvable groups are amenable,
\item{iii)} Free discrete groups (with more than one generator)
are not amenable,
\item{iv)}
 If there is an exact sequence of topological groups
  $0\to A \to B \to C\to 0$, then $B$ is amenable iff $A$ and $C$
  are amenable.
\item{v)} Closed subgroups of amenable groups are amenable.

Using these facts we can determine completely (for $char(k)=0$) when
$G(k)$ is amenable.

\Lemma
Let $k$ be a local field, $G$ a $k$-group, $R$ the radical of $G$.
For $char(k)=0$ the follwing properties are equivalent:
\item{1)}
$G/R$ is $k$-anisotropic, \ie $G/R$ does not contain any positive-dimensional
$k$-split torus;
\item{2)} $G/R(k)$ is compact;
\item{3)} $G(k)/R(k)$ is compact;
\item{4)} $G(k)$ is amenable.
\endproclaim
\Proof
For $1)\!\iff\! 2)$ see \cite M. The equivalence
$2)\!\iff\! 3)$ follows from Theorem C,
because $G(k)/R(k)$ is closed in $G/R(k)$ (\cite {BT,3.18.}).
The implication $3)\implies 4)$ is a direct consequence of the above
listed properties on amenable groups (i), ii) \& iv)).
Finally, if $G/R(k)$ is not compact, then there exists a discrete subgroup
of $G(k)$ containing a free subgroup (Th. A \& B).
But an amenable group cannot contain a closed discrete non-commutative
free subgroup. Hence $4)$ implies $2)$.
\qed

There is a similar result for
arbitrary connected real Lie groups, see \cite{Z, 4.1.9 on p.62}.

In positive characteristic the radical $R$ is not necessarily defined
over $k$, hence $1)$ and $2)$ do not make sense.
Even if $R$ is defined over $k$, it is not clear whether
$G(k)/R(k)$ has finite index in $G/R(k)$. Therefore it is not clear
whether $2)$ and $3)$ are equivalent.
But at least the implication $3)\implies 4)$
is true in any characteristic.

\Proposition
Let $G$ be a Zariski-connected
 $k$-group defined over a local field $k$.
Assume that $G(k)$ is amenable and that there exists a Zariski-dense
discrete subgroup $\Gamma$ of $G(k)$.

Then $G$ is solvable.
\endproclaim
\Proof
Amenability of $G$ implies amenability of $\Gamma$.
By the result of Tits either $\Gamma$ contains a free subgroup or
a solvable subgroup of finite index.
Since the former is ruled out by amenability of $\Gamma$,
$\Gamma$ contains a solvable subgroup $\Gamma_0$ of finite index which
is still Zariski-dense in $G$. Hence $G$ is solvable.
\qed
\Corollary
Let $G$ be a Zariski-connected
$k$-group defined over a local field $k$, $R$ its radical.
Assume that $G/R(k)$ is compact, but $G\ne R$.

Then $G(k)$ does not contain any Zariski-dense discrete subgroup.
\endproclaim

\section{An obstruction for solvable groups}
We will now derive an obstruction to the existence of discrete
Zariski-dense subgroups in solvable groups.

For $A,B\subset G$ let $[A,B]$ denote the subgroup generated by
the commutators $aba^{-1}b^{-1}$.
\Lemma
Let $A,B$ algebraic subgroups in $G$ and $\Gamma$ a subgroup of $G$ such
that $A\cap\Gamma$ and $B\cap\Gamma$ are Zariski-dense in $A$ resp. $B$.

Then $[A,B]\cap\Gamma$ is Zariski-dense in $[A,B]$.
\endproclaim
\Proof
Let $C$ denote the set of commutators $aba^{-1}b^{-1}$ and
$C^n=C\cdot\ldots\cdot C$. For $n$ sufficiently large the natural
morphism $(A\times B)^n\to C^n\to [A,B]$
is dominant and therefore maps the Zariski-dense subset
$\((A\cap\Gamma)\times(B\cap\Gamma)\)^n$ onto a Zariski-dense subset
$\Lambda$ in $[A,B]$.
\qed
\def\C{{\cal C}}
\Definition
For a group $G$ let $\C(G)$ denote the smallest
collection of subgroups of $G$ such that
$G\in\C(G)$ and $[A,B]\in\C(G)$ for all $A,B\in\C(G)$.
\endproclaim
The collected $\C(G)$ contains in particular all subgroups of the
derived and the central series.
Every $H\in\C(G)$ is a normal subgroup and for a $k$-group $G$
every $H\in\C(G)$ is defined over $k$ (\cite{B,p.58}).
For a Zariski-connected group $G$ every $H\in\C(G)$ is again
Zariski-connected.

We will use this notation to deduce
an obstruction to the existence of discrete Zariski-dense
subgroups.
\Lemma
Let $H$ be a Zariski-connected one-dimensional $k$-group, $k$ a local field,
$\Gamma\subset H(k)$ an infinite discrete subgroup.
Let $A$ denote the group of all $k$-group automorphisms of $H$
stabilizing $\Gamma$.

Then $A$ is finite.
\endproclaim
\Proof
If $H$ is $\bar k$-isomorphic to the multiplicative group $G_m$, then
the group of {\sl all} automorphisms of $H$ is finite
(because $z\to z^{-1}$ is the only non-trivial automorphism of $G_m$).
Hence we may assume that $H\simeq G_a$. This isomorphism
is given over some finite extension field $k'$ of $k$.
(\cite{B,Th.10.9 \& Remark}).
We may replace $k$ by $k'$, \ie we may assume that $H\simeq G_a$
as $k$-groups.
Now any $k$-group automorphism of $H$ is given by $\mu_\lambda:
z\mapsto\lambda z$
for some $\lambda\in k^*$.
Since $k$ is locally compact and $\Gamma$ discrete, there exists
an element $\gamma_0\in\Gamma$ such that $\abs{\gamma_0}\le\abs g$
for all $g\in\Gamma\setminus\{0\}$.
This easily implies that $\abs\lambda=1$ for all $\mu_\lambda\in A$.
Now the quotient map from $A$ to the $A$-orbit $A(\gamma_0)$
is injective and $A(\gamma_0)\subset\{x:\abs{x}=\abs{\gamma_0}\}$.
Since $\{x:\abs{x}=c\}$ is compact for all $c$, it follows that
$A$ is finite.
\qed

\Proposition
Let $k$ be a local field, $G$ a Zariski-connected $k$-group, $H\in\C(G)$.
Assume that $H$ is one-dimensional and not central.
Then $G$ does not admit any discrete Zariski-dense subgroup.
\endproclaim
\Remark
(1) Though not stated explicitly, the assumptions of this proposition imply
that $G$ is solvable. To see this, note that for each $H\in\C(G)$ there
exists a subgroup $I$ in the derived series of $G$ such that $I\subset H$.
It follows that $G$ is solvable as soon as there exists a solvable
subgroup $H\in\C(G)$.

(2) As we will see below, solvable non-commutative p-adic groups never
admit Zariski-dense discrete subgroups.
Hence this proposition is interesting only for $k=\s R$, $k=\s C$ and
$char(k)>0$.
\endremark
\Proof of Proposition 3
Assume that there exists a discrete Zariski-dense subgroup
$\Gamma$ of $G(k)$.
Thanks to Lemma 2 it is clear that $A\cap\Gamma$ must be Zariski-dense
in $A$ for all $A\in\C(G)$.
Hence $H\cap\Gamma$ is infinite.
Now $\Gamma$ acts by conjugation on $H$, stabilizing $H\cap\Gamma$.
With the help of the preceding lemma this implies that $\Gamma$
contains a subgroup $\Gamma_0$ of finite index such that $\Gamma_0$
centralizes $H$.
Since $H$ is {\sl not} central in $G$, this contradicts the assumption
that $\Gamma$ is Zariski-dense in $G$.
\qed
We will now apply this to Borel groups.
\Corollary
Let $S$ be a Zariski-connected simple $k$-group
and assume that there exists a Borel
subgroup $G$ defined over $k$.
Then $G(k)$ does not admit any discrete Zariski-dense subgroup.
\endproclaim
\Proof
Let $H$ be the one-dimensional unipotent subgroup of $G$ corresponding
to the maximal root.
It is easy to check that $H$ is not central and $H\in\C$.
\qed

As we will see below, simplicity is an essential condition for this result.
We will demonstrate that a Borel group $B$ of
the semisimple complex linear-algebraic group $SL_2(\s C)\times SL_2(\s C)$
{\sl does} admit discrete Zariski-dense subgroups.
\bigskip
\vfill
\eject
\section{Unipotent Groups in characteristic zero}
\Lemma
Let $U$ be a unipotent $k$-group, $char(k)=0$.
Then every (non-trivial) element of $U$ is of infinite order
and $U(k)$ contains a Zariski-dense subgroup generated by $r=dim(U)$
elements.
\endproclaim
\Proof
Recall that $U$ is $k$-split, \ie admits a sequence of $k$-subgroups
$G=G_0\supset\ldots\supset G_s=\{e\}$ with $\dim G_i/G_{i+1}=1$.
\cite{B, Cor. 15.5 (ii) on p.205}.
Furthermore for $char(k)=0$ every one-dimensional unipotent $k$-group
is $k$-isomorphic to $G_a$ (\cite{B,Th.10.9 \& Remark below}).
Using these facts the assertions of the lemma follows easily by induction
on $dim(U)$.
\qed
\Lemma
Let $U$ be a commutative unipotent $k$-group, $char(k)=0$,
$g\in U\setminus\{e\}$ and
$\s Z(g)\subset U$ the subgroup generated by $g$.
Then the Zariski-closure of $\s Z(g)$ is one-dimensional.
\endproclaim

\Proof
This follows immediately from the fact that $U$ is
$\bar k$-isomorphic to a vector group $G_a^n$
(\cite{S,p.171}).
\qed
By induction one obtains the following consequence.
\Corollary
Let $U$ be a commutative unipotent $k$-group, $char(k)=0$,
$\Gamma$ a Zariski-dense subgroup.
Then $\Gamma$ can not be generated
by less then $dim(U)$ elements.
\endproclaim

Combining these two lemmata we obtain
\Proposition
Let $U$ be a commutative unipotent $k$-group, $char(k)=0$.
Then $U$ is $k$-isomorphic to the vector group $(G_a)^r$.
\endproclaim
\Proof
Let $r=dim(U)$ and $g_1,\ldots,g_r\in U(k)$ generators of a Zariski-dense
subgroup. Let $\s Z(g_i)$ denote the subgroups generated by each $g_i$ and
$H_i$ the Zariski-closure of the $\s Z(g_i)$.
Since $g_i\in U(k)$, the $H_i$ are defined over $k$
(\cite{B,AG 14.4}).
Furthermore the $H_i$ are one-dimensional and $k$-isomorphic to $G_a$.
Since $U$ is commutative, the $k$-subgroups $H_i$
induce a $k$-group homomorphism $\rho$ from $H=\Pi H_i\simeq(G_a)^r$
to $U$. Now $\rho$ maps a Zariski-dense subgroup of $H(k)$
onto the Zariski-dense
subgroup $\Gamma\subset U(k)$ generated by the $g_i$ and is therefore
dominant and defined over $k$ (\cite{BT,1.4.}). 
A dominant $k$-group morphism between $k$-groups of the same dimension
is surjective with finite kernel.
But $char(k)=0$ implies that $(G_a)^r$ has no non-trivial finite subgroup.
Thus $\rho$ is bijective.
Again using $char(k)=0$ it follows
that $\rho$ is an isomorphism.
Finally, the inverse $\rho^{-1}$ maps the Zariski-dense set $\Gamma
\subset U(k)$ into $H(k)$ and is therefore likewise defined over $k$.
\qed
Such a statement does not hold in positive characteristic,
consider
\eg  Witt groups \cite{S,VII} or Examples 3 and 4 in Section 9.
\Corollary   1
Let $U$ be a commutative unipotent $k$-group, $k$ an archimedean local
field (\ie $\s R$ or $\s C$). Then $U$ admits a Zariski-dense
discrete subgroup.
\endproclaim

\Corollary  2
Let $U$ be a commutative unipotent $k$-group, $k$ a non-archimedean
local field of characteristic zero, $\Gamma$ a finitely generated
subgroup of $U(k)$.
Then $\Gamma$ is relatively compact in $U(k)$.
\endproclaim
\Proof
This follows from the ultrametric condition via $U(k)\simeq G_a(k)^r$.
\qed

\Corollary 3
Let $G$ be a commutative unipotent $k$-group, $k$ non-archimedean local field
of
characteristic zero.

Then $G(k)$ contains no discrete subgroups except $\{e\}$.
\endproclaim
\Proof
Using induction on $dim(G)$ it is easy to prove that
for $char(k)=0$ a unipotent
$k$-group cannot contain a non-trivial finite subgroup.
Therefore the preceding corollary implies that $\{e\}$ is the only
finitely generated discrete subgroup of $G(k)$.
Finally note that any discrete group must contain a finitely generated
subgroup.
\qed
\Corollary 4
Let $G$ be a unipotent $k$-group,
$k$ a non-archimedean local field of characteristic zero.

Then $G(k)$ contains no discrete subgroups except $\{e\}$.
\endproclaim
\Proof
Assume the contrary. Since every element in $G(k)$ is of infinite order,
such a discrete subgroup would contain a subgroup $\Gamma$
isomorphic to $\s Z$. But then the Zariski-closure of $\Gamma$ would be
a commutative unipotent $k$-group, thus contradicting the preceding
corollary.
\qed

\Proposition
Let $G$ be a Zariski-connected solvable $k$-group for a non-archimedean
local field $k$ with $char(k)=0$.
Then any discrete subgroup $\Gamma\subset G(k)$ is commutative.
In particular $G(k)$ cannot admit any discrete Zariski-dense subgroup
unless $G$ is commutative.
\endproclaim
\Proof
By standard results on solvable groups the commutator group $\cal D(G)$
of $G$ is unipotent. Thus for any discrete subgroup $\Gamma\subset G(k)$
we have $\Gamma\cap\cal D(G)(k)=\{e\}$, hence $\Gamma$ is commutative.
This implies that the Zariski-closure of $\Gamma$ is commutative, too.
\qed

Now we turn to unipotent groups defined over archimedean fields.
Here the question of the existence of discrete Zariski-dense subgroups
has been settled by the following result of Malcev.
\Theorem D (Malcev, \cite{Mc}, see also \cite R)
Let $G$ be a unipotent $k$-group, $k=\s R$.
Then $G(k)$ admits a Zariski-dense discrete subgroup if and only
if $G$ may be defined over $\s Q$.
Each such discrete Zariski-dense subgroup is cocompact.
\endproclaim
In Malcev's article the condition that $G$ may be defined
over $\s Q$ is replaced by the property that the Lie algebra may be defined
over $\s Q$. However, for unipotent groups in characteristic zero the
exponential map gives an isomorphism (as $k$-varieties) of the group
and its Lie algebra. Using the Campbell-Hausdorff formula it follows
that $G$ may be defined over $\s Q$ if and only if the Lie algebra can be
defined over $\s Q$.

Malcevs result immediately implies the following criterion for complex
unipotent groups.
\Corollary
Let $G$ be a unipotent $k$-group, $k=\s C$.
By "restriction of scalars" $G(\s C)$ is isomorphic
(as topological group) to $\tilde G(\s R)$
for a unipotent group $\tilde G$ defined over $\s R$.
Fix a continuous group isomorphism $\phi:\tilde G(\s R)\to G(\s C)$.

$G(\s C)$ admits a Zariski-dense discrete subgroup if and only if
$\tilde G$ admits a unipotent subgroup $H$ defined over $\s Q$
such that $\phi(H(\s R))$ is Zariski-dense in $G$.
\endproclaim
For the convenience of the reader we reformulate this in the language
of Lie groups.
\Corollary
Let $G$ be a unipotent complex Lie group.
Then $G$ admits a discrete subgroup which is dense in the algebraic
Zariski-topology if and only if
$G$ contains a real Lie subgroup $H$ such that
\item{1)}
The structure constants for $\Lie(H)$ are rational numbers for a suitable
base,
\item{2)}
$\Lie(G)$ is the smallest complex vectorsubspace of $\Lie(G)$
containing $\Lie(H)$.
\endproclaim
\section{Real solvable groups}
We will now use our results on real unipotent groups in order to deduce
a statement about arbitrary real solvable groups.
\Proposition
Let $G$ be a Zariski-connected solvable $k$-group, $k=\s R$.
The following conditions are necessary
for the existence of a
Zariski-dense discrete subgroup.
\item{(i)}
The commutator group $\cal D(G)$ is defined over $\s Q$.
\item{(ii)}
$G$ is unimodular.
\endproclaim

\Proof
Let $\Gamma$ be a discrete Zariski-dense subgroup in $G(\s R)$.
Then the commutator group
$\cal D(\Gamma)$ is Zariski-dense in $\cal D(G)$.
Hence property $(i)$ is necessary.
Moreover this implies that $\cal D(\Gamma)$ is cocompact in $\cal D(G)(\s R)$.
This ensures unimodularity of $G$.
\qed
These conditions are not sufficient.
For instance, let
$$G(\s R)=\left\{\pmatrix{
\lambda^{-4} & w &&& \cr
&1&&& \cr
&&\lambda^2 & x & z \cr
&&&\lambda & y \cr
&&&& 1 \cr}:\lambda\in\s R^*;x,y,z\in\s R\right\}$$
$G$ fulfills all the conditions of the theorem, but likewise fulfills
the obstruction criterion deduced in Proposition 3:
$G''$ is one-dimensional, but not central.
\section{Unipotent groups in positive characteristics}
For local fields in positive characteristics
we have only fragmentary results.
A local field in positive characteristics is isomorphic to $\s F_q((t))$
for some $q=p^n$. Such a field contains infinite discrete subrings,
\eg the ring generated by the elements $t^k$ for $k\le 0$.
Hence it is easy to give examples of unipotent $k$-groups which do
admit discrete Zariski-dense subgroups.
However these discrete Zariski-dense subgroups are never finitely generated:
\Lemma
Let $k$ be a field of positive characteristic, $U$ a
Zariski-connected unipotent $k$-group
and $\Gamma$ a finitely generated subgroup of $U(k)$.

Then $\Gamma$ is finite.
\endproclaim
\Proof
There is no loss in generality, if we assume $k$ to be algebraically
closed.
Then the lemma follows easily by induction, since any one-dimensional
unipotent $k$-group is isomorphic to $G_a$.
(In addition,
one has to use the fact that for a finitely generated {\sl nilpotent}
group every subgroup is again finitely generated
\cite{R,Th.2.7.}).
\qed

This is in strong contrast to the situation in characteristic zero.
For a local field $k$ with $char(k)=0$, a
discrete subgroup of a solvable $k$-group is always finitely
generated (This follows from \cite{R, Prop.3.8}).

There also exist unipotent commutative $k$-groups in positive characteristic
which do not contain any Zariski-dense discrete subgroup.
\Example 3
Let $k=\s F_p((t))$ and $G$ the one-dimensional unipotent $k$-group
defined by $G=\{(x,y)\in G_a\times G_a:x^p-x=ty^p\}$.%
\footnote*{This group has been studied by M. Rosenlicht
\cite{Ro,p.46} for a different purpose.}
Then $G(k)$ is compact and therefore does not admit any infinite discrete
subgroup.
To check compactness, let $x=\sum_i a_it^i$, $y=\sum_i b_it^i$.
An explicit calculation shows that
$(x,y)\in G$ if and only if $a_k^p=a_{kp}$ and $-a_{kp+1}=b_k^p$
for all $k\in\s Z$ and $a_k=0$ for all $k$ with
$k \hbox{ \sl mod }p\not\in\{0,1\}$.
This implies $a_k=0=b_k$ for all $k<0$. Hence $G(k)$ is a closed subgroup
of $\cal O\times\cal O$, where $\cal O=\{x:\abs x\le 1\}$ denotes
the (compact) additive group of local integers.
Thus $G(k)$ is compact.
\endproclaim
Even for non-compact unipotent groups it is possible that there exist
no discrete Zariski-dense subgroup.
\Example 4
Let $k=\s F_p((t))$.
Let $U=G_a\times W_2$ where $W_2$ is the two-dimensional Witt group,
\ie $W_2$ is $\s A^2$ as variety with the group multiplication given
by $(x,y)\cdot(z,w)=(x+z,y+w+F(x,z))$ for
$F(x,z)={1\over p}\left(x^p+y^p-(x+y)^p\right)$. (See \cite{S,VII.2} for
more about Witt groups).
As a variety $U=\s A^3$. Now let $H=\{(x,y,z)\in U:x^p-x=ty\}$.
As a $k$-group $H$ is an extension of $G_a$ by the group $G$ studied
in the previous example.
Now $H(k)$ is non-compact, because it contains $G_a(k)$.
Nevertheless there are no discrete Zariski-dense subgroups in $H(k)$.
To see this, consider the group morphism $\phi:x\mapsto x^p$, where
$x^p$ denotes the $p$-th power with respect to group multiplication
in $H$. By standard results on Witt groups $\phi$ is a dominant morphism
from $H$ to $A=\{(0,0,z)\}\subset U$ with $A\subset \ker\phi$.
Now $H(k)/A(k)$ is compact, hence $\phi(H(k))$ is compact.
Since $\phi(\Gamma)\subset\Gamma$, it follows that $\phi(\Gamma)$ is finite
for every discrete subgroup $\Gamma\subset H(k)$.

But $\phi(\Gamma)$ must be Zariski-dense in $\phi(H)=A$ for every
Zariski-dense subgroup $\Gamma\subset H$.
\endproclaim

\section{General commutative groups}
We will now deal with commutative groups which are not necessarily
unipotent.
The following decomposition theorem is a centerpiece for this
investigation.
\Theorem D
Let $G$ be a Zariski-connected commutative $k$-group, $char(k)=0$.
Then
$G$ admits a decomposition
$G=G_u\times G_c\cdot G_i$
with $G_u$ unipotent, $G_c$ $k$-anisotropic torus
and $G_i$ a $k$-split torus (\ie
$G_i$ is $k$-isomorphic to $G_m(k)^r$).
$G_c\cap G_i$ is finite and $G_u\cap(G_c\cdot G_i)=\{e\}$.
All the groups $G_u$, $G_c$ and $G_i$ are defined over $k$.
Moreover this decomposition is functorial, \ie any group
morphism $\rho:G\to H$ between commutative $k$-groups will
map $G_u$, $G_i$ and $G_c$ into $H_u$, $H_i$ resp. $H_c$.
\endproclaim
The same statement holds for positive characteristic except that
$G_u$ is only $k$-closed and not necessarily defined over $k$.
\Proof
See \cite{B,p.121ff and p.137ff}.
\qed
\Proposition
Let $G$ be a Zariski-connected
commutative $k$-group, where $k$ is a non-archi\-me\-dean
local field with $char(k)=0$.

Then $G(k)$ admits a discrete Zariski-dense subgroup
if and only if
\item{(1)} $G\ne G_c$ and
\item{(2)} $\dim(G_i)\ge\dim(G_u)$.
\endproclaim
\Proof
$G_c(k)$ is the maximal Zariski-connected compact subgroup of $G(k)$, hence
condition (i) is clearly necessary.
Now let $\Gamma$ be a finitely generated discrete subgroup of $G(k)$.
Let $a_i$ denote the generators and $\pi_u:G\to G_u$ the natural projection.
By Cor.2 to Prop.\ 4 the image of $\Gamma$ in $G_u$ under the natural
projection
is relatively-compact.
Thus the projection of $\Gamma$ to $G_c\times G_u$ has a relatively compact
image. It follows that the natural projection $\pi_i:G\to G_i$ has finite
kernel if restricted to $\Gamma$.
Now there is a proper continuous group homomorphism
$G_i(k)\to\s Z^r$ with $r=\dim(G_i)$ induced by the logarithm of the
absolute value on the local field $k$.
Hence there is a group homomorphism $\Gamma\to\s Z^r$ with finite kernel.
It follows that any discrete subgroup $\Lambda$ of $G(k)$ is finitely generated
with $rank(\Lambda)\le\dim(G_i)$.
Due to the corollary to Lemma 5
a finitely generated Zariski-dense subgroup $\Lambda_u$ of $G_u$
must fulfill $rank(\Lambda_u)\ge\dim(G_u)=n$.
Therefore condition $(ii)$ is necessary for
the existence of discrete Zariski-dense subgroups.
On the other hand by the functoriality of the decomposition a subgroup
$\Gamma$ is Zariski-dense in $G$ if and only if all the projections to $G_c$,
$G_i$ and $G_u$ are Zariski-dense.
Now $G_c(k)$ admits an element generating a Zariski-dense subgroup
(\cite {B,Prop.8.8 \& Remark below}) and $G_u(k)$ admits
a Zariski-dense subgroup generated by $dim(G_u)=n$ elements
(Lemma 4).
For any $x\in k^*$ with $\abs{x}\ne 1$ we obtain a Zariski-dense
discrete subgroup $\Gamma$ of $G_i\simeq(k^*,\cdot)^m)$ by
$\Gamma=\{(x^{n_1},\ldots,x^{n_m}):n_i\in\s Z\}$.
Thus $dim(G_i)\ge\max(dim(G_i),1)$ is sufficient for the existence of
discrete Zariski-dense subgroups in $G(k)$.
\qed
\Remark
For local fields of positive characteristic the "if"-part of the
proposition still holds. However, the "only if"-part breaks down,
\eg $G_a(k)$ admits a Zariski-dense discrete subgroup.
\endremark
The following is well-known and easy to prove.
\Lemma
Let $G$ be a Zariski-connected commutative $k$-group for $k=\s R$ or $k=\s C$.
For $k=\s R$ the group
$G(k)$ admits a discrete cocompact Zariski-dense subgroup
if and only if $G(k)$ is non-compact.
For $k=\s C$ the group $G(k)$ is never compact and always admits a
discrete Zariski-dense discrete cocompact subgroup.
\endproclaim

\section{Metabelian groups}
We start with an auxiliary remark.
\Lemma
Let $\rho:\s C^*\to GL(V)$ be a rational representation on a complex
vector space $V$.
Assume that all the weights are non-zero and let $\lambda\in\s C^*$.
Then there exists a number $n\in\s N$ such that either all
or none of the eigenvalues
of $\rho(\lambda)$ are real.
\endproclaim
\Proof
The eigenvalues of $\rho(\lambda)$ are $\lambda^{k_1},\ldots\lambda^{k_n}$
for some $k_1,\ldots k_n\in\s Z\setminus\{0\}$.
If $\lambda^{k_i}\in\s R$ for an $i$,
then $\rho(\lambda^{k_i})$ and $\rho(\lambda^{-k_i})$ have
only real eigenvalues.
\qed
We use this to derive a necessary conition for the existence of Zariski-dense
discrete subgroups in certain metabelian groups.
\Proposition
Let $\rho:\s C^*\to GL(V)$ be a rational representation for which all the
weights are distinct and non-zero.
Let $G=\s C^*\sdr_\rho V$ be the induced semidirect product.
Assume that $G$ admits a discrete Zariski-dense subgroup $\Gamma$.

Then $G$ is unimodular, \ie $\rho(\s C^*)\subset SL(V)$.
\endproclaim
\Proof
Consider the projection $\tau:G\to\s C^*\simeq G/G'$.
Now $\tau(\Gamma)$ is Zariski-dense in $\s C^*$.
We claim that $\tau(\Gamma)$ is not relatively compact in $\s C^*$.
Assume the contrary. Then $S=\overline{\tau(\Gamma)}$ is compact.
Thus the action of $\Gamma$ by conjugation on $V$ factors through a
compact group action on $V$. Hence the orbits of $\Gamma$ acting by
conjugation on $\Gamma\cap V$ are relatively compact (and discrete),
hence finite. Since $\Lambda=\Gamma\cap V$ is a finitely generated
group, it follows that $\Gamma$ contains a subgroup of finite index
which centralizes $\Lambda$. This is impossible, because $\Lambda$
is Zariski-dense in $V=[G,G]$ and $V$ is not central in $G$.

Thus $\tau(\Gamma)\subset\s C^*$ is not relatively compact
and in particular contains
an element which is not of finite order, \ie not a root of unity.
Let $\lambda$ be such an element.
By the preceding lemma we may assume either all or none of the eigenvalues
of $\rho(\lambda)$ are real.
Now let $V=\oplus_\omega V_\omega$ be a decomposition of $V$ into eigenspaces
of $\rho(\lambda)$.
By assumption they are one-dimensional.
Let $\pi_\omega:V\to V_\omega$ denote the respective projections.
Since $\Gamma\cap V$ is Zariski-dense in $V$, it contains an element
$v$ such that $\pi_\omega(v)\ne 0$ for all $\omega$.
Let $\Sigma$ denote the smallest $\rho(\lambda)$-invariant subgroup of
$V$ containing $v$. Then $W=\Sigma\tensor_{\s Z}\s R$ is a real subvectorspace
of $V$, which is again invariant.
It follows that $W=\oplus_\omega(W\cap
V_\omega)$. If the eigenvalues are complex, then $W=V$.
Since $\Sigma$ is a lattice in $W$, in this case $\Sigma$ is a lattice in $V$.
If all eigenvaluies are real, then $W$ is a totally real subvectorspace
with $V=W\oplus iW$. Thus in this case $\Sigma\oplus i\Sigma$ is a lattice
in $V$.
In any case $V$ admits a lattice, which is stable under $\rho(\lambda)$.
Therefore $\rho(\lambda)\in SL(V)$. Since the group generated by $\lambda$
is Zariski-dense in $\s C^*$, it follows that $\rho(\s C^*)\subset SL(V)$.
Thus $G$ is unimodular.
\qed
However, even in this special
case unimodularity is not sufficient for the existence of
a discrete Zariski-dense subgroup.
\Example 5
Let $G=\s C^*\sdr_\rho\s C^3$ with the weights of $\rho$ given
by $(2,-1,-1)$, \ie $\rho(\lambda)(x_1,x_2,x_3)=(\lambda^2x_1,
\lambda^{-1}x_2,\lambda^{-1}x_3)$.
Let $V=\s C^3=V_2\oplus V_{-1}$ where $V_\alpha$ denotes the weight
space for $\alpha$.
Assume that $G$ admits a Zariski-dense discrete subgroup $\Gamma$.
Then $\Gamma'$ is Zariski-dense in $G'=V$.
Hence $\Gamma$ contains an element $\gamma$ which is contained in $V$ but
neither in $V_{2}$ nor in $V_{-1}$.
Then $\gamma=\gamma_2+\gamma_{-1}$ with $0\ne\gamma_\alpha\in V_\alpha$.
The elements $\gamma_\alpha$ span a two-dimensional subvectorspace $W$ of $V$.
Consider the natural projection $\tau:G\to G/G'\simeq\s C^*$.
By the considerations in the proof of the above proposition there exists
an element $\delta\in\Gamma$ such that $\abs{\tau(\delta)}>1$
and $\tau(\delta)\in\s R$ or $(\tau\delta)^2\not\in\s R$
Let $\Lambda$ denote the subgroup of $\Gamma$ generated by $\gamma$ and
$\delta$. Then by the same reasoning
as in the above proof $\Lambda'$ or $\Lambda'+i\Lambda'$ must be lattice
in $W$ which is stable under $\rho(\delta)$.
But this is impossible, because $\abs{\det\left(\rho(\delta)|_W\right)}>1$.
\endproclaim

\section{A number-theoretical construction}
A series of metabelian linear-algebraic groups over $\s C$ with
discrete cocompact subgroups may be constructed by the following
number-theoretic approach which generalizes a construction of
Otte and Potters \cite{OP}.
\footnote*{Otte and Potters studied the special case where $K$ is a totally
real numberfield. For this case our theorem is equivalent to
\cite{OP,3.2.}.}

\Theorem
Let $K$ be a number field, $\cal O$ the ring of algebraic integers,
$\cal O^*$ the multiplicative group of units, $r_1$ the number of
real imbeddings, $r_2$ the number of pairs of conjugate complex
embeddings and $r=r_1+r_2-1$.

Then there is a semidirect product
$G=T\sdr_\rho V$ of a torus $T=(\s C^*)^m$ with
and a vector group $V=\s C^{r+1}$ with
$m\in\{r,r+1\}$ and $\{e\}=\ker\rho:T\to GL(V)$, such that
$G$ contains a Zariski-dense discrete subgroup $\Gamma$ isomorphic
to $\cal O^*\sdr\cal O$.

For $r_1=r_2=1$, we obtain $m=2$ and $G$ is isomorphic to a Borel group
in $SL_2(\s C)\times SL_2(\s C)$.

For totally real $K$ there exists a discrete cocompact subgroup
$\Gamma_0$ isomorphic to $\cal O^*\sdr(\cal O\times\cal O)$.
Furthermore for totally real $K$, $\rho(T)$ is a maximal torus
in $SL(V)$.
\endproclaim
\Proof
Let $\zeta_1,\ldots,z_{r_1}$ and $\xi_1,\bar\xi_1,\ldots
\xi_{r_2},\bar\xi_{r_2}$ denote the real resp. pairs of conjugate
complex imbeddings of $K$ in $\s C$.
Then $\phi=(\zeta_1,\ldots,\zeta_{r_1},\xi_1,\ldots,\xi_{r_2})$ embed
$K$ into the complex vector space $V=\s C^{r+1}$
such that $W=K\tensor_{\s Q}\s R$
is embedded into $V$ as a real subvectorspace with $\phi(W)+i\phi(W)=V$.
Furthermore $\phi(\cal O)$ is a cocompact lattice in $\phi(W)$ and
therefore a discrete Zariski-dense subgroup of $V$.
The action of $K^*$ on $K$ by multiplication induces an action on $V$.
For each $x\in K^*$ this action is a diagonalizable endomorphism of $V$
with eigen-values $\zeta_1(x),\ldots,\zeta_{r_1}(x),\xi_1(x),\ldots
\xi_{r_2}(x)$. It follows that $K^*$ acts on $V$ as a subgroup of
a maximal torus $T$ of $GL_r(\s C)$. Thus we obtain an injective
group homomorphism $\tau:K^*\to T$. Now $\tau(\cal O^*)$ stabilizes the
discrete Zariski-dense subgroup $\cal O$ of $V$.
It is easy to see that this implies that $\tau(\cal O^*)$ is discrete in $T$.
Hence we obtain a discrete group $\Gamma=\tau(\cal O^*)\sdr\phi(\cal O)$
which is Zariski-dense in its Zariski-closure $G=\bar\Gamma$.
Clearly $V\subset G\subset T\sdr V$ and $G'=[G,G]=V$.
By the theorem of Dirichlet $\cal O^*$ is isomorphic to
a direct product of a finite
abelian group $A$ and $\s Z^r$.
Hence $\dimc(G/V)\ge r$.
On the other hand $\dimc(G/V)\le r+1$, since $\dimc(T)=\dimc(V)=r+1$.

Now let us consider the special case $r_1=r_2=1$.
Then for $\alpha\in\cal O^*$ we obtain
$\abs{\det(\tau(\alpha))}=\abs{\zeta_1(\alpha)}\abs{\xi_1(\alpha)}$.
Now $\abs{N_{K,\s Q}(\alpha)}=1$ implies that either $G$ is not unimodular or
$\abs{\xi_1(\alpha)}=1=\abs{\zeta_1(\alpha)}$.
However, the second alternative would imply that $\tau(\cal O^*)$ is
relatively compact in $GL(V)$, which is impossible (see the proof of
Prop. 9). Hence $G$ is not unimodular.
By Prop. 9 it follows that $\tau(G)\subset GL(V)$ can not be one-dimensional.
Hence $\tau(G)$ must be a maximal torus in $GL(V)=GL_2(\s C)$, which
implies that $G$ is a Borel group in $S=SL_2(\s C)\times SL_2(\s C)$.

Finally let us consider the special where $K$ is totally real, \ie
$r_2=0$
(This is the case studied in
\cite{OP}).
Then $\phi(W)$ is totally real in $V$ and we obtain a cocompact
lattice in $V$ by $\Lambda=\phi(\cal O)+i\phi(\cal O)$.
Furthermore in this case $\det(\tau(x))=
=N_{K,\s Q}(x)$ for all $x\in K^*$.
Since the norm $N_{K,\s Q}(x)$ equals $1$ or $-1$ for all units $x\in\cal
O^*$, it follows that $\tau(\cal O^*)$ admits a subgroup $\Delta$ of index $2$
or $1$ which is contained in $SL(V)$.
Now $\Delta\sdr\Lambda$ is a discrete cocompact subgroup in $T_0\sdr V$
where $T_0=T\cap SL(V)$.
\qed

\References
\item{[B]} Borel, A.: Linear algebraic groups.
Second enlarged edition. Springer 1991.

\item{BS} Borel, A.; Serre, J.P.:
Th\'eor\`emes de finitude en cohomologie galoisienne.
\sl Comm. Math. Helv. \bf 39\rm, 111--164 (1964)

\item{BT} Borel, A.; Tits, J.:
Homomorphismes "abstraits" de groupes alg\'ebriques simples.
Ann. Math. \bf 97\rm, 499-571 (1973)

\item{G} Greenleaf, F.:
Invariant means on topological groups.
Van Nostrand, New York 1969

\item{[Mc]} Malcev, A.:
On a class of homogeneous spaces.
\sl Izvestiya Akad. Nauk SSSR Ser. Math. \bf 13 \rm (1949)/
\sl AMS Transl. no. \bf 39 \rm (1951)

\item{M} Margulis, G.A.:
Discrete Subgroups of Semisimple Lie Groups.
Springer Berlin Heidelberg New York 1989

\item{OP} Otte, M.; Potters, J.:
Beispiele homogener Mannigfaltigkeiten.
\sl Manu. math. \bf 10\rm, 117--127 (1973)

\item{[R]} Raghunathan, M.S.:
Discrete subgroups of Lie groups.
\sl\rm Erg. Math. Grenzgeb. \bf 68 \rm, Springer (1972) \*

\item{Ro} Rosenlicht, M.:
Some rationality questions on algebraic groups.
\sl Annali di Math. \bf 43\rm, 25-50 (1957)

\item{S} Serre, J.P.:
Algebraic Groups and Class Fields.
Springer 1988.

\item{T} Tits, J.:
Free subgroups in linear groups.
\sl J. Algebra \bf 20\rm, 250--270 (1972)

\item{Z} Zimmer, R.J.:
Ergodic Theory and Semisimple Groups.
Birkh\"auser 1984

\endarticle
\endinput